\documentclass{amsart}
\usepackage{amssymb}
\usepackage{amsfonts}
\usepackage{amsmath}
\usepackage{graphicx}
\usepackage[pagebackref]{hyperref}
\theoremstyle{plain}
\newtheorem{thm}{Theorem}[section]
\newtheorem{corollary}[thm]{Corollary}

\theoremstyle{definition}

\theoremstyle{remark}
\newtheorem{remark}[thm]{Remark}
\newtheorem{example}[thm]{Example}

\numberwithin{equation}{section}
\def\1{{\rm (1)}}
\def\2{{\rm (2)}}
\def\3{{\rm (3)}}
\def\4{{\rm (4)}}
\def\5{{\rm (5)}}
\begin{document}

%%%%%%%%%%%%%%%%%%%%%%%%%%%%%%%%%%%%%%%%%%%%%%%%%%%%%%%%%
%%%%%%%%%%%%%%%%%%%%%%%%%%%%%%%%%%%%%%%%%%%%%%%%%%%%%%%%%
\title[Well-centered Overrings of a Commutative Ring in Pullbacks and Trivial Extensions]
{Well-centered Overrings of a Commutative Ring in Pullbacks and Trivial Extensions}

%%%%%%%%%%%%%%%%%%%%%%%%%%%%%%%%%%%%%%%%%%%%%%%%%%%%%%%%%
%%%%%%%%%%%%%%%%%%%%%%%%%%%%%%%%%%%%%%%%%%%%%%%%%%%%%%%%%
\author{N. Mahdou}

\address{Department of Mathematics, Faculty of Sciences and Technology of Fez, Box 2202, University S. M. Ben Abdellah, Fez,
Maroc}

\email{mahdou@hotmail.com}

\author{A. Mimouni}

\address{Department of Mathematics and Statistics,
King Fahd University of Petroleum \& Minerals, P. O. Box 278,
Dhahran 31261, Saudi Arabia}

\email{amimouni@kfupm.edu.sa}

\thanks{The second named author is supported by KFUPM}


%%%%%%%%%%%%%%%%%%%%%%%%%%%%%%%%%%%%%%%%%%%%%%%%%%%%%%%%%%%%%%%%%%%%%%%%%%%%%%%%%%%%%%%%%%%%%%%%%%%%%%%%%
\begin{abstract} Let $R$ be a commutative ring with identity and $T(R)$ its total quotient ring.
We extend the notion of well-centered overring of an integral domain to an arbitrary commutative ring and
we investigate the transfer of this property to different extensions of commutative rings in both integral
and non-integral cases. Namely in pullbacks and trivial extensions. Our
aim is to provide new classes of commutative rings satisfying this property and to shed light on some open questions
raised by Heinzer and Roitman in \cite{HR}.

\end{abstract}

\maketitle
%%%%%%%%%%%%%%%%%%%%%%%%%%%%%%%%%%%%%%%%%%%%%%%%%%%%%%%%%%%%%%%%%%%%%%%%%%%%%%%%%%%%%%%%%%%%%%%%%%%%%%%%
%%%%%%%%%%%%%%%%%%%%%%%%%%%%%%%%%%%%%%%%%%%%%%%%%%%%%%%%%%%%%%%%%%%%%%%%%%%%%%%%%%%%%%%%%%%%%%%%%%%%%%%%
\section{Introduction}

Throughout, $R$ is a commutative ring with identity and $Tot(R)$
denotes its total quotient ring (in case where $R$ is an integral
domain, we denote by $L$ its quotient field). A regular element is
an element which is not a zerodivisor and a regular ideal is an ideal
containing a regular element. By an integral ideal $I$, we mean an ideal
$I$ such that $I\subseteq R$ and by fractional ideal, we mean a nonzero
$R$-submodule $E$ of $Tot(R)$ such that $dE\subseteq R$ for a regular
element $d$ of $R$. Finally by overring, we mean a ring $T$ such that $R\subseteq T\subseteq Tot(R)$.\\

In 2004, Heinzer and Roitman \cite{HR} introduced and studied the notion of well-centered overrings of an integral domain, that is, an overring $T$ of a domain $R$ is well-centered on $R$ if each element of $T$ is an associate in $T$ of an element of $R$, equivalently, if each principal ideal of $T$ is generated by an element of $R$ (note that a similar definition in case of a valuation overring of $R$ is given in \cite{Gr}). Since a localization of $R$ is both flat over $R$ and well-centered on $R$, the authors studied when the converse holds, i.e., when a flat and well-centered overring of $R$ is a localization of $R$. They proved that in general the converse does not hold, however in some particular cases such as simple extensions (i.e. $T=R[b]$) or if $T$ is finitely generated on $R$ and $R$ is Noetherian, the answer is yes. The following three questions was raised:\\

\noindent{\bf Question 1} (\cite[Question 4.7]{HR}). Under what conditions on an integral domain $R$ is every finitely generated well-centered
overring of $R$ a localization of $R$?\\

\noindent{\bf Question 2} (\cite[Question 4.11]{HR}). Under what conditions on an integral domain $R$ is every flat overring of $R$ well-centered on $R$?\\

\noindent{\bf Question 3} (\cite[page 454]{HR}). Whether a finitely generated flat well-centered
overring of an integral domain $R$ is always a localization of $R$?\\

The purpose of this paper is to continue the investigation of the notion of well-centered overrings. We first extend this notion to an arbitrary commutative ring. In Section 2, we deal with well-centered overrings of an integral domain issued from pullbacks. Our motivation is an example constructed by Heinzer and Roitman in order to produce a simple extension of a Noetherian domain $R$ (so well-centered on $R$) but which is not flat over $R$. We prove that for a pullback $R$ issued from a diagram of type $(\square)$, if $T$ is local, then it is well centered over $R$, and if each intermediate overring $S$ between $R$ and $T$ is well-centered on $R$, then $K=T/M$ is algebraic over $D=R/M$. Moreover, if $T$ is local and $D=k$ is a field, then every intermediate ring between $R$ and $T$ is well-centered on $R$ if and only if $K$ is algebraic over $k$ (Theorem~\ref{Pull.1}). This lead us to construct a families of integral domains $R$ such that every overring of $R$ is well-centered on $R$ in both Noetherian and non-Noetherian cases. Also we prove that in the context of $PVDs$ (pseudo-valuation domains), a well-centered overring of a Mori domain is a Mori domain (Corollary~\ref{Pull.3}).\\

In Section 3 we investigate the transfer of ``well-centered property" to trivial extensions. Our aim is to construct a families of overrings of a commutative ring $R$ (which is not necessarily a domain) such that every finitely generated well-centered overring of $R$ is a localization of $R$ and every flat overring of $R$ is well-centered on $R$, but for which $R$ is not even a coherent ring. This shows that a characterization of commutative rings satisfying the three pre-cited questions is so far from being realized. Our main results state that if $A\subseteq B$ is an extension of rings,  $E$ is a $B$-module, $T: =B \propto E$ is the trivial ring extension of $B$ by $E$, and $R =A \propto E$ is the trivial ring extension of $A$ by $E$, then $T$ is well-centered on $R$ if and only if $B$ is well-centered on $A$ (Theorem~\ref{WCTE.2}); and if $A$ is a ring, $T(A)$ its total quotient ring, $E$ a $T(A)$-module and $R :=A \propto E$, then every finitely generated well-centered overring of $R$ is a localization of $R$ if and only if every finitely generated well-centered overring of $A$ is a localization of $A$, and every flat overring of $R$ is well-centered on $R$ if and only if every flat overring of $A$ is well-centered on $A$ (Theorem~\ref{WCTE.4}).\\
\bigskip
%%%%%%%%%%%%%%%%%%%%%%%%%%%%%%%%%%%%%%%%%%%%%%%%%%%%%%%%%%%%%%%%%%%%
%%%%%%%%%%%%%%%%%%%%%%%%%%%%%%%%%%%%%%%%%%%%%%%%%%%%%%%%%%%%%%%%%%%%
%%%%%%%%%%%%%%%%%%%%%%%%%%%%%%%%%%%%%%%%%%%%%%%%%%%%%%%%%%%%%%%%%%%%

\noindent{\bf Definition 1.1}(\cite[page 435]{HR})
Let $R$ be an integral domain and $T$ an overring of $R$.
We say that $T$ is well-centered on $R$ if for each $b\in T$, there exists a unit $u\in T$ such that $ub=a\in R$.
Thus $T$ is well-centered on $R$ iff each element of $T$ is an associate in $T$ of an element of $R$ iff each principal ideal of $T$ is
generated by an element of $R$.\\

We extend Heinzer-Roitman's definition to an arbitrary commutative ring in the following way:\\

\noindent{\bf Definition 1.2} Let $R$ be a commutative ring with identity and $T$ an overring of $R$. We say that
$T$ is well-centered on $R$ if for each $b\in T$, there exists a unit $u\in T$ such that $ub=a\in R$.
\bigskip
%%%%%%%%%%%%%%%%%%%%%%%%%%%%%%%%%%%%%%%%%%%%%%%%%%%%%%%%%%%%%%%%%%%%%%%%%%%%%%%%%%%%%%%%%%%%%%%%%%%%%%%%%%%%%%%%%%%%%
%%%%%%%%%%%%%%%%%%%%%%%%%%%%%%%%%%%%%%%%%%%%%%%%%%%%%%%%%%%%%%%%%%%%%%%%%%%%%%%%%%%%%%%%%%%%%%%%%%%%%%%%%%%%%%%%%%%%%
%%%%%%%%%%%%%%%%%%%%%%%%%%%%%%%%%%%%%%%%%%%%%%%%%%%%%%%%%%%%%%%%%%%%%%%%%%%%%%%
%%%%%%%%%%%%%%%%%%%%%%%%%%%%%%%%%%%%%%%%%%%%%%%%%%%%%%%%%%%%%%%%%%%%%%%%%%%%%%%
\section{Pullbacks}\label{Pull}

The purpose of this section is to investigate the notion of well-centered overrings
in pullback constructions. Our work is motivated by a result due to Heinzer and Roitman
\cite[Proposition 3.9]{HR}, where they proved that if $T$ is a domain of the form $T=K+M$, where
$K$ is a subfield of $T$, $M$ a maximal ideal of $T$ and $D$ a domain contained in $K$, and $R=D+M$, then
$T$ is always well-centered on $R$.
This construction is used to prove the existence of a simple well-centered proper overring $B$ of a Noetherian domain $A$
that is a sublocalization of $A$, but which is not flat (\cite[Examples 3.10 and 3.11]{HR}). Our aim is to generate new
families of well-centered overrings. First, let us fix the notation for the rest of
this section and recall some useful properties of pullbacks.\\
Let $T$ be an integral domain, $M$ a maximal ideal of $T$, $K=T/M$
its residue field and $D$ a subring of $K$. Let $R$ be the pullback
of the following diagram

\[\begin{array}{ccl}
R            & \longrightarrow                 & D\\
\downarrow   &                                 & \downarrow\\
T            & \stackrel{\phi}\longrightarrow  & K=T/M
\end{array}\]

We assume that $R\subsetneq T$ and we refer to this diagram as a
diagram of type $(\square)$ and if $qf(D)=K$, we refer to the
diagram as a diagram of type $(\square^{*})$. The case where $T=V$
is a valuation domain is of crucial interest, we shall
refer to this as a classical diagram.\\
Recall that $(R:T)=M$ is a prime ideal of $R$, $R/M\simeq D$ and $R$ and $T$ have
the same quotient field. Moreover, if $T$ is local, then every ideal of $R$ is comparable
(under inclusion) to $M$, and $R$ is local if and only if $T$ and $D$ are local.
For more details on general pullbacks, we refer the reader to \cite{F, FG, GH}
and \cite{BG} for classical ``D+M" constructions.\\

\begin{thm}\label{Pull.1} Let $R$ be the pullback of the diagram of type $(\square)$. Then\\
\1 If $T$ is local, then $T$ is well-centered on $R$.\\
\2 If each intermediate overring $S$ between $R$ and $T$ is well-centered on $R$, then $K$ is algebraic over $D$.
Moreover, if $T$ is local and $D=k$ is a field, then every intermediate ring between $R$ and $T$ is well-centered on $R$ if and only if
$K$ is algebraic over $k$.
\end{thm}

\begin{proof} \1 Let $b\in T$. If $b\in M$, we are done. If $b\not\in M$, then $b$ is
a unit of $T$ and in this case just take $a=1$, so that $bT=aT=T$.\\
\2 Assume that each overring $S$ of $R$ such that $R\subseteq S\subseteq T$ is well-centered on $R$.
Let $\lambda\in K\setminus D$ and let $x\in T$ such that $\phi(x)=\lambda$. Set $S=\phi^{-1}(k[\lambda])$, where $k=qf(D)$.
Then $S$ is well-centered on $R$ and so there exists $y\in R$ such that $xS=yS$. Since $\lambda\not =0$, $x\not\in M$ and so $y\not\in M$.
But $xy^{-1}\in U(S)$ implies that $xy^{-1}z=1$ for some $z\in S$, and so $xz=y$. Hence $\lambda\phi(z)=\phi(y)\in D\setminus \{0\}$. So $(\lambda\phi(z))^{-1}\in k$ and therefore $\lambda^{-1}=(\lambda\phi(z))^{-1}\phi(z)\in k[\lambda]$. Hence $\lambda$ is algebraic over $k$, as desired.\\
Finally assume that $T$ is local, $D=k$ is a field and $K$ is algebraic over $k$. Then every intermediate ring $S$ between $R$ and $T$ is of the form $S=\phi^{-1}(F)$ for some intermediate field $F$ between $k$ and $K$. But since $T$ is local, then so is $S$ and by part \1, $S$ is well-centered on $R$.
\end{proof}
%%%%%%%%%%%%%%%%%%%%%%%%%%%%%%%%%%%%%%%%%%%%%%%%%%%%%%%%%%%%%%%%%%%%%%%%%%%%%%%%%%%%%%%%%%%%%%%%%%%%%%%%%%%%%%%%%%%%%%
%%%%%%%%%%%%%%%%%%%%%%%%%%%%%%%%%%%%%%%%%%%%%%%%%%%%%%%%%%%%%%%%%%%%%%%%%%%%%%%%%%%%%%%%%%%%%%%%%%%%%%%%%%%%%%%%%%%%%%
The following example shows that if $T$ is not local or if $T$ is not of the form $T=K+M$,
then $T$ is not necessarily well-centered on $R$.
\bigskip
%%%%%%%%%%%%%%%%%%%%%%%%%%%%%%%%%%%%%%%%%%%%%%%%%%%%%%%%%%%%%%%%%%%%%%%%%%%%%%%%%%%%%%%%%%%%%%%%%%%%%%%%%%%%%%%%%%%%%%
%%%%%%%%%%%%%%%%%%%%%%%%%%%%%%%%%%%%%%%%%%%%%%%%%%%%%%%%%%%%%%%%%%%%%%%%%%%%%%%%%%%%%%%%%%%%%%%%%%%%%%%%%%%%%%%%%%%%%%
\begin{example}\label{Pull.2} Let $\mathbb{Q}$ be the field of rational numbers and $X$ an indeterminate over
$\mathbb{Q}$. Set $T=\mathbb{Q}[X]$, $M=(X^{2}+2)T$ and $K=T/M=\mathbb{Q}(\sqrt{2})$, and
let $R$ be the pullback arising from the following diagram of canonical maps:\\
\[\begin{array}{ccl}
R            & \longrightarrow                 & \mathbb{Q}\\
\downarrow   &                                 & \downarrow\\
T=\mathbb{Q}[X]            & \stackrel{\phi}\longrightarrow  & K=T/M=\mathbb{Q}(\sqrt{2})
\end{array}\]
We claim that $T$ is not well-centered on $R$. Otherwise, there exists $f\in R$ such that $fT=XT$. Hence
$Xf^{-1}$ is a unit of $T$ and so $Xf^{-1}=a$ for some nonzero element $a\in \mathbb{Q}$. Therefore $f=aX$ and so
$a\sqrt{2}=\phi(aX)=\phi(f)\in \mathbb{Q}$, a contradiction.
\end{example}
%%%%%%%%%%%%%%%%%%%%%%%%%%%%%%%%%%%%%%%%%%%%%%%%%%%%%%%%%%%%%%%%%%%%%%%%%%%%%%%%%%%%%%%%%%%%%%%%%%%%%%%%%%%%%%%%%%%%%%%%%%%%
%%%%%%%%%%%%%%%%%%%%%%%%%%%%%%%%%%%%%%%%%%%%%%%%%%%%%%%%%%%%%%%%%%%%%%%%%%%%%%%%%%%%%%%%%%%%%%%%%%%%%%%%%%%%%%%%%%%%%%%%%%%%
In \cite[Example 3.24]{HR}, Heinzer and Roitman showed that a well-centered overring of a Mori domain is not, in general, a Mori domain. Our next result shows that in the context of $PVDs$ (pseudo-valuation domains) the Mori property is preserved by well-centered property.\\
\bigskip
%%%%%%%%%%%%%%%%%%%%%%%%%%%%%%%%%%%%%%%%%%%%%%%%%%%%%%%%%%%%%%%%%%%%%%%%%%%%%%%%%%%%%%%%%%%%%%%%%%%%%%%%%%%%%%%%%%%%%%%%%%%%%%%%%%%%%%%%%%%%%%%
%%%%%%%%%%%%%%%%%%%%%%%%%%%%%%%%%%%%%%%%%%%%%%%%%%%%%%%%%%%%%%%%%%%%%%%%%%%%%%%%%%%%%%%%%%%%%%%%%%%%%%%%%%%%%%%%%%%%%%%%%%%%%%%%%%%%%%%%%%%%%%%

\begin{corollary}\label{Pull.3} Let $R$ be a $PVD$, $V$ its associated valuation overring and $M$ its maximal ideal.\\
\1 Each overring of $R$ is well-centered on $R$ if and only if $V/M$ is  algebraic over $R/M$.\\
\2 If $R$ is Mori, then each well-centered overring of $R$ is Mori.
\end{corollary}

\begin{proof}
Proposition 2.6 of \cite{AD} characterizes PVDs in terms of
pullbacks. The aforementioned proposition states that $R$ is a PVD
if and only if $R=\phi^{-1}(k)$ for some subfield $k$ of $K=V/M$,
where $V$ is the associated valuation overring of $R$, $M$ its
maximal ideal and $\phi$ the canonical homomorphism from $V$ onto
$K$.\\
\1 Follows immediately from Theorem~\ref{Pull.1} and the fact that each overring $S$ of $R$ containing $V$ is of the form $S=V_{P}$ for some prime ideal $P$ of $R$ (since each overring of $R$ is comparable to $V$, \cite[Theorem 2. 1]{BG}).\\
\2 Assume that $R$ is a Mori domain and let $T$ be a well-centered overring of $R$. Since $V$ is a $DVR$ and $T$ is comparable to $V$, we may assume that $T\subsetneq V$. Then $T=\phi^{-1}(B)$ for some integral domain $B$ satisfying $k\subseteq B\subseteq K$. In virtue of \cite[Theorem 4.18]{GH}, it suffices to show that $B$ is a filed. Let $\lambda\in B\setminus \{0\}$ and let $x\in T$ such that $\phi(x)=\lambda$. Since $T$ is well-centered on $R$, there exists $y\in R$ such that $xT=yT$. Since $\lambda\not =0$, $x\not\in M$ and so $y\not\in M$.
But $xy^{-1}\in U(T)$ implies that $xy^{-1}z=1$ for some $z\in T$, and so $xz=y$. Hence $\lambda\phi(z)=\phi(y)\in k\setminus \{0\}$. So $(\lambda\phi(z))^{-1}\in k$ and therefore $\lambda^{-1}=(\lambda\phi(z))^{-1}\phi(z)\in B$. Hence $B$ is a field, as desired.
\end{proof}
\bigskip
%%%%%%%%%%%%%%%%%%%%%%%%%%%%%%%%%%%%%%%%%%%%%%%%%%%%%%%%%%%%%%%%%%%%%%%%%%%%%%%%%%%%%%%%%%%%%%%%%%%%%%%%%%%%%%%%%%%%%%%%%%%%%%%%%%%%%%%%%%%%%%%%%%%%%%%
%%%%%%%%%%%%%%%%%%%%%%%%%%%%%%%%%%%%%%%%%%%%%%%%%%%%%%%%%%%%%%%%%%%%%%%%%%%%%%%%%%%%%%%%%%%%%%%%%%%%%%%%%%%%%%%%%%%%%%%%%%%%%%%%%%%%%%%%%%%%%%%%%%%%%%%
\begin{remark}\label{Pull.4} Let $R$ be a $PVD$, $V$ its associated valuation overring, $M$ its maximal ideal, $K=V/M$, $k=R/M$ and suppose that $K$ is algebraic over $k$. By Corollary~\ref{Pull.3}, each overring of $R$ is well-centered on $R$. The Noetherianity of $R$ depends now on $V$, that is, if $V$ is a $DVR$, then $R$ is Noetherian and if $V$ is not a $DVR$, then $R$ is not Notherian. This leads us to construct a families of integral domains $R$ such that every overring of $R$ is well-centered on $R$ in both Noetherian and non-Noetherian contexts. Moreover, for such classes, no proper overring of $R$ can be a localization of $R$.
\end{remark}
\bigskip

%%%%%%%%%%%%%%%%%%%%%%%%%%%%%%%%%%%%%%%%%%%%%%%%%%%%%%%%%%%%%%%%%%%%%
%%%%%%%%%%%%%%%%%%%%%%%%%%%%%%%%%%%%%%%%%%%%%%%%%%%%%%%%%%%%%%%%%%%%%
\section{Well-centered property in trivial ring extensions}\label{WCTE}

Let $A$ be a  commutative ring and let $U(A)$ denote the set of units of $A$.
Let $E$ be an $A$-module and  $R =A \propto  E$, the
set of pairs $(a,e)$ with $a\in A$ and $e\in E$, under
coordinatewise addition and under an adjusted multiplication
defined by $(a,e)(a',e')=(aa',ae'+a'e),$  for all $a,a'\in A,
e,e'\in E$. Then $R$ is called the trivial ring extension of
$A$ by $E$. We recall that the units of $R$ are of the form $(a,e)$, where $a$ is unit of $A$.\\
Trivial ring extensions have been studied extensively; the work is
summarized in Glaz \cite{G1} and Huckaba \cite{H}. These extensions have
been useful for solving many open problems and conjectures in both
commutative and non-commutative ring theory, see for instance \cite{G1}, \cite{H}, \cite{KM}.\\
In this section we investigate the transfer
of well-centered property between a ring and its trivial ring extensions. For this, we start by giving a complete description
of the form of overrings of a ring issued by trivial extension.\\
\bigskip
%%%%%%%%%%%%%%%%%%%%%%%%%%%%%%%%%%%%%%%%%%%%%%%%%%%%%%%%%%%%%%%%%%%%%%%%%%%%%%%%%%%%%%%%%%%%%%%%%%%%%%%%%%%%%%%%%%%%%%%%%%%%%%
%%%%%%%%%%%%%%%%%%%%%%%%%%%%%%%%%%%%%%%%%%%%%%%%%%%%%%%%%%%%%%%%%%%%%%%%%%%%%%%%%%%%%%%%%%%%%%%%%%%%%%%%%%%%%%%%%%%%%%%%%%%%%%
\begin{thm}\label{WCTE.1}
Let $A$ be a ring, $T(A)$ be a total ring of quotient of $A$, $E$ be a $T(A)$-module,
$R :=A \propto E$, and let $S$ be a ring. Then $S$ is an overring of $R$ if and only if
there exists an overring $B$ of $A$ such that $S =B \propto E$.
\end{thm}

\begin{proof} Let $B$ be an overring of $A$ and set $S :=B \propto E$. Our aim is to show that $S$
is an overring of $R$, that is, $S \subseteq Tot(R)$. \\
Let $(b,e) \in S$ where $b \in B$ and $e \in E$. Hence, $b =a/s$ where $a \in A$ and $s$ is a regular element of $A$ since $B \subseteq Tot(A)$. Therefore, $(s,0)(b,e) =(sb,se) =(a,se) \in R$ and it remains to show that $(s,0)$ is a regular element of $R$. \\
Let $(c,d) \in R$ such that $(s,0)(c,d) =(0,0)$, which means that $sc =0$ in $A$ and $sd =0$
in $E$. Hence, $c =0$ since $s$ is a regular element of $A$. On the other hand,
$sd =0$ implie that $d =0$ since $s$ is an invertible element of $Tot(A)$ and $E$ is a $Tot(A)$-module,
as desired. \\

Conversely, assume that $S$ is an overring of $R$, that is, $R \subseteq S \subseteq T(R)$. Our aim is to show that
$T(R) =T(A) \propto E$ and this suffices to show that $S$ has the form $S =B \propto E$, where $B$ is
an overring of $A$. \\
Let $U$ be a set of regular element of $A$. Then it is a multiplicative set of both $A$ and $R$; and we have
$U^{-1}R =(U^{-1}A) \propto (U^{-1}E) =T(A) \propto E$ since $E$ is a $T(A)$-module and $U^{-1}T(A) =T(A)$. But,
we can easily show that every element of $U$ is a regular element of $R$. Therefore,
$(R \subseteq ) U^{-1}R =T(A) \propto E \subseteq T(R)$ and so it remains to show that $T(A) \propto E$ is a total ring,
that is, every element of $T(A) \propto E$ is invertible or zero-divisors. \\
Let $(a,e)$ be a regular element of  $T(A) \propto E$. Then $a \not= 0$ (since if $a =0$, then $(0,e)(0,f) =0_{R}$ for each $f \in E$ which means
that $(0,e)$ is not a regular element of  $T(A) \propto E$). Now, it is easy to show that $a$ is a regular element of $T(A)$ since $(a,e)$ is a regular element of $T(A) \propto E$ and so $a$ is invertible in $T(A)$
(since $T(A)$ is a total ring of quotient). Therefore,
$(a,e)$ is invertible in $T(A) \propto E$ and this completes the proof of  the Theorem.
\end{proof}
\bigskip
%%%%%%%%%%%%%%%%%%%%%%%%%%%%%%%%%%%%%%%%%%%%%%%%%%%%%%%%%%%%%%%%%%%%%%%%%%%%%%%%%%%%%%%
%%%%%%%%%%%%%%%%%%%%%%%%%%%%%%%%%%%%%%%%%%%%%%%%%%%%%%%%%%%%%%%%%%%%%%%%%%%%%%%%%%%%%%%

Now, we give our first main result in this section.\\
\bigskip
%%%%%%%%%%%%%%%%%%%%%%%%%%%%%%%%%%%%%%%%%%%%%%%%%%%%%%%%%%%%%%%%%%%%%%%%%%%%%%%%%%%%%%%%%%%%%%%%%%%%%%%%%%%%%%%%%%%%%%%%%%%%%%%%%%%%%%%%%%%%%%%%%%%%
%%%%%%%%%%%%%%%%%%%%%%%%%%%%%%%%%%%%%%%%%%%%%%%%%%%%%%%%%%%%%%%%%%%%%%%%%%%%%%%%%%%%%%%%%%%%%%%%%%%%%%%%%%%%%%%%%%%%%%%%%%%%%%%%%%%%%%%%%%%%%%%%%%%%
\begin{thm}\label{WCTE.2}
Let $A \subseteq B$ be an extension of rings, $E$ be an $B$-module, $T: =B \propto E$ be the trivial ring extension of
$B$ by $E$, and let $R =A \propto E$ be the trivial ring extension of $A$ by $E$. Then
$T$ is well-centered on $R$ if and only if $B$ is well-centered on $A$.
\end{thm}

\begin{proof} Assume that $T$ is well-centered on $R$ and let $c \in B$. Hence,
there exists a unit element $(b,e)$ in $T$, where $b \in B$ and $e \in E$, such that
$(c,0)(b,e) (=(bc, ce)) \in R$. Thus, $b$ is a unit element in $B$ (since $(b,e)$ is a unit in $T$)
and $bc \in A$; this means that $B$ is well-centered on $A$. \\

Conversely, assume that $B$ is well-centered on $A$ and let $(b,e) \in T$,
where $b \in B$ and $e \in E$. If $b =0$, then $(0,e) \in R$ and the result is clear. Assume that $b \not= 0$.
Since $B$ is well-centered on $A$ and $b \in B$, there exists $u \in U(B)$ such that $bu \in A$.
Hence $(b,e)(u,0) =(ub, ue) \in R$ and $(u,0) \in U(T)$. This means that $T$ is well-centered on $R$ and this completes the proof.
\end{proof}
\bigskip
%%%%%%%%%%%%%%%%%%%%%%%%%%%%%%%%%%%%%%%%%%%%%%%%%%%%%%%%%%%%%%%%%%%%%%%%%%%%%%%%%%%%%%%%%%%%%%%%%%%%%%%%%%%%%%%%%%%%%%%%%%%%%%%%%%%%%%%%%%%%%%%%%%%
%%%%%%%%%%%%%%%%%%%%%%%%%%%%%%%%%%%%%%%%%%%%%%%%%%%%%%%%%%%%%%%%%%%%%%%%%%%%%%%%%%%%%%%%%%%%%%%%%%%%%%%%%%%%%%%%%%%%%%%%%%%%%%%%%%%%%%%%%%%%%%%%%%%
The following example show that the condition ``$E$ is a $T(A)$-module" in Theorem~\ref{WCTE.1} and
the condition ``$T: =B \propto E$ and $R: =A \propto E$" in Theorem~\ref{WCTE.2} are necessary.\\
\bigskip
%%%%%%%%%%%%%%%%%%%%%%%%%%%%%%%%%%%%%%%%%%%%%%%%%%%%%%%%%%%%%%%%%%%%%%%%%%%%%%%%%%%%%%%%%%%%%%%%%%%%%%%%%%%%%%%%%%%%%%%%%%%%%%%%%%%%%%%%%%%%%%
%%%%%%%%%%%%%%%%%%%%%%%%%%%%%%%%%%%%%%%%%%%%%%%%%%%%%%%%%%%%%%%%%%%%%%%%%%%%%%%%%%%%%%%%%%%%%%%%%%%%%%%%%%%%%%%%%%%%%%%%%%%%%%%%%%%%%%%%%%%%%%
\begin{example}\label{WCTE.3}
Let $A$ be a ring, $B$ and $C$ two overrings of $A$ such that $A \subseteq B \subsetneqq C\subseteq T(A)$.
Set $R :=A \propto B$ and  $S :=B \propto C$. Then: \\
{\bf 1)} $S$ is an overring of $R$. \\
{\bf 2)} $S$ is never well-centered on $R$ (even if $A$ is a valuation domain).
\end{example}

\begin{proof} {\bf 1)} It is clear that $S$ is an overring of $R$ since $T(R) =T(A) \propto T(A)$. \\
{\bf 2)} We claim that $S$ is never well-centered on $R$. Deny, let $(0,c) \in S$,
where $c \in C - B$. Then there exists $(u,v) \in U(B \propto C) (=U(B) \propto C)$
such that $(u,v)(0,c) (=(0,uc)) \in R$ since $S$ is well-centered on $R$. Hence, $uc \in B$
and so $c \in B$ since $u \in U(B)$, a contradiction. Therefore,
$S$ is never well-centered on $R$.
\end{proof}
\bigskip
%%%%%%%%%%%%%%%%%%%%%%%%%%%%%%%%%%%%%%%%%%%%%%%%%%%%%%%%%%%%%%%%%%%%%%%%%%%%%%%%%%%%%%%
%%%%%%%%%%%%%%%%%%%%%%%%%%%%%%%%%%%%%%%%%%%%%%%%%%%%%%%%%%%%%%%%%%%%%%%%%%%%%%%%%%%%%%%%%%%%%%%%%%
In \cite[Theorem 4.5]{HR} Heinzer and Roitman proved that for a Pr\"ufer domain $R$ with Noetherian
spectrum, every finitely generated well-centered overring of $R$ is a localization of $R$, and  for a Noetherian
domain $R$, every finitely generated flat overring of $R$ is a localization of $R$. They raised the
following open questions:\\

\noindent{\bf Question 1} (\cite[Question 4.7]{HR}). Under what conditions on an integral domain $R$ is every finitely generated well-centered
overring of $R$ a localization of $R$?\\

\noindent{\bf Question 2} (\cite[Question 4.11]{HR}). Under what conditions on an integral domain $R$ is every flat overring of $R$ well-centered on $R$?\\

\noindent{\bf Question 3} (\cite[page 454]{HR}). Whether a finitely generated flat well-centered
overring of an integral domain $R$ is always a localization of $R$?\\

In what follows, and by extending the above questions to an arbitrary commutative ring (that is not necessarily a domain), we study the transfer of these questions in trivial ring extensions. Further, we construct new classes of arbitrary commutative rings $R$ such that every finitely generated well-centered overring and every flat overring is a localization of $R$. Our motivation is based on a simple remark on valuation domains
as a simple and well-known class of integral domains satisfying the above questions.\\
\bigskip
%%%%%%%%%%%%%%%%%%%%%%%%%%%%%%%%%%%%%%%%%%%%%%%%%%%%%%%%%%%%%%%%%%%%%%%%%%%%%%%%%%%%%%%%%%%%%%%%%%%%%%%%%%%%%%%%%%%%%%%
%%%%%%%%%%%%%%%%%%%%%%%%%%%%%%%%%%%%%%%%%%%%%%%%%%%%%%%%%%%%%%%%%%%%%%%%%%%%%%%%%%%%%%%%%%%%%%%%%%%%%%%%%%%%%%%%%%%%%%%
\begin{thm}\label{WCTE.4}
Let $A$ be a ring, $T(A)$ its total quotient ring, $E$ a $T(A)$-module and let
$R :=A \propto E$. Then:\\
\1 Every finitely generated well-centered overring of $R$ is a localization of $R$ if and only if every finitely generated well-centered overring of $A$ is a localization of $A$.\\
\2 Every flat overring of $R$ is well-centered on $R$ if and only if every flat overring of $A$ is well-centered on $A$.\\
\3 Every finitely generated flat well-centered overring of $R$ is a localization of $R$ if and only if every finitely generated flat well-centered overring of $A$ is a localization of $A$.
\end{thm}

\begin{proof} \1 Assume that every finitely generated well-centered overring of $A$ is a localization of $A$ and let $S$ be a finitely generated well-centered overring of $R$. By Theorem~\ref{WCTE.1}, $S$ has the form: $S =B \propto E$, where $B$ is an overring of $A$.
Also, $B$ is well-centered on $A$ by Theorem~\ref{WCTE.2} since $S$ is well-centered on $R$. Trivially $B$ is a finitely generated $A$-module (since if $S =\displaystyle\sum _{i=1}^{n}R(a_{i}, e_{i})$,
where $n$ is a positive integer and $(a_{i}, e_{i}) \in S$ for each $i =1, \ldots ,n$, then $B =\displaystyle\sum _{i=1}^{n}Aa_{i}$). Therefore, $B$ is a localization of $A$. Let $U$ be a multiplicative set
of $A$ such that $B =U^{-1}A$. Hence, $U^{-1}R =U^{-1}(A \propto E) =(U^{-1}A) \propto (U^{-1}E)
=B \propto E =S$ (since $U^{-1}E =S^{-1}BE =BE =E$), as desired. \\

Conversely, assume that every finitely generated well-centered overring of $R$ is a localization of $R$ and let $B$ be a finitely generated well-centered overring of $A$. By Theorem~\ref{WCTE.1}, $S =B \propto E$ is an overring of $R$.
Also, $S$ is well-centered on $R$ by Theorem~\ref{WCTE.2} since $B$ is well-centered on $A$. On the other hand, $S$ is  a finitely generated $R$-module. Indeed, set $B =\displaystyle\sum _{i=1}^{n}Aa_{i}$,
where $n$ is a positive integer and $a_{i} \in B$ for each $i =1, \ldots ,n$. Then
$\displaystyle\sum _{i=1}^{n}R(a_{i},0) =(\displaystyle\sum _{i=1}^{n}Aa_{i}) \propto (\displaystyle\sum _{i=1}^{n}Aa_{i}E) =B \propto (BE) =B \propto E=S$ since $E$ is a $B$-module (as $E$ is a $T(A)$-module
and $T(A) =T(B)$). Therefore, $S$ is a finitely generated well-centered overring of $R$ and so
$S$ is a localization of $R$. Let $U$ be a multiplicative set
of $R$ such that $S =U^{-1}R$. Then $U \subseteq U(B \propto E) =U(B) \propto E$. \\
Let $U_{0} :=\{s \in A | (s,v) \in U$ for some $v \in E\}$. It is easy to see that
$U_{0}$ is a multiplicative set of $A$ and $U_{0}^{-1}A =B$ and this completes the proof of \1.\\

\2 Assume that every flat overring of $A$ is well-centered on $A$ and let $S$ be a flat
overring of $R$. By Theorem~\ref{WCTE.1}, $S$ is of the form $S =B \propto E$, where $B$ is an overring of $A$.
But $S$ is a flat $A$-module since $S$ is a flat $R$-module and $R$ is a flat $A$-module.
Therefore, $B$ is a flat $A$-module since $S \cong B \times E$ as $A$-modules. Hence
$B$ is well-centered on $A$, and so $S$ is well-centered on $R$ by Theorem~\ref{WCTE.2}, as desired. \\

Conversely, let $B$ be a flat overring of $A$. By Theorem~\ref{WCTE.1}, $S =B \propto E$ is an overring of $R$. We claim that
$S$ is $R$-flat. Indeed, $B \otimes_{A} R$ is a flat $R$-module since $B$ is a flat $A$-module. But,
$B \otimes_{A} R =B \propto E =S$ since $R \cong A \times E$ as $A$-modules and
$B \otimes_{A} E =BE =E$ (since $E$ is a $T(A)$-module and $T(A) =T(B)$). Therefore, $S$ is a flat overring of $R$ and so $S$ is well-centered on $R$. Hence $B$ is well-centered on $A$ by Theorem~\ref{WCTE.2} as desired.\\

\3 Similar to \1 and \2.
\end{proof}
\bigskip
%%%%%%%%%%%%%%%%%%%%%%%%%%%%%%%%%%%%%%%%%%%%%%%%%%%%%%%%%%%%%%%%%%%%%%%%%%%%%%%%%%%%%%%%%%%%%%%%%%%%%%%%%%%%%%%%%%%%%%%%%%%%%%%%%%%%%%%
%%%%%%%%%%%%%%%%%%%%%%%%%%%%%%%%%%%%%%%%%%%%%%%%%%%%%%%%%%%%%%%%%%%%%%%%%%%%%%%%%%%%%%%%%%%%%%%%%%%%%%%%%%%%%%%%%%%%%%%%%%%%%%%%%%%%%%%
Now, we are ready to construct a new families of arbitrary commutative rings satisfying the above questions and that are not even coherent.
This shows that suitable characterizations for commutative rings satisfying one of the above questions are so far from being realized.\\
%%%%%%%%%%%%%%%%%%%%%%%%%%%%%%%%%%%%%%%%%%%%%%%%%%%%%%%%%%%%%%%%%%%%%%%%%%%%%%%%%%%%%%%%%%%%%%%%%%%%%%%%%%%%%%%%%%%%%%%%%%%%%%%%%%%%%%%
%%%%%%%%%%%%%%%%%%%%%%%%%%%%%%%%%%%%%%%%%%%%%%%%%%%%%%%%%%%%%%%%%%%%%%%%%%%%%%%%%%%%%%%%%%%%%%%%%%%%%%%%%%%%%%%%%%%%%%%%%%%%%%%%%%%%%%%
\bigskip

\begin{example}
Let $V$ be a valuation domain which is not a field, $L :=qf(V)$, $E$ be an $L$-vector space and let $R :=V \propto E$. Then: \\
{\bf 1)} Every finitely generated well-centered overring of $R$ a localization of $R$.\\
{\bf 2)} Every flat overring of $R$ is well-centered on $R$. \\
{\bf 3)} $R$ is not coherent. \\
\end{example}

\begin{proof} {\bf 1)} and  {\bf 2)} follow immediately from Theorem~\ref{WCTE.4}.\\
{\bf 3)} Let $e (\not= 0) \in E$ and $(0,e) \in R$. Then
$(0:(0,e)) =0 \propto E$ is not a finitely generated
ideal of R (since $E$ is a $K$-vector space and $K$ is not a finitely generated $V$-module).
Therefore, $R$ is not a coherent ring, as asserted.
\end{proof}
%%%%%%%%%%%%%%%%%%%%%%%%%%%%%%%%%%%%%%%%%%%%%%%%%%%%%%%%%%%%%%%%%%%%%%%%%%%%%%%%%%%%%%%
%%%%%%%%%%%%%%%%%%%%%%%%%%%%%%%%%%%%%%%%%%%%%%%%%%%%%%%%%%%%%%%%%%%%%%%%%%%%%%%%%%%%%
%%%%%%%%%%%%%%%%%%%%%%%%%%%%%%%%%%%%%%%%%%%%%%%%%%%%%%%%%%%%%
%%%%%%%%%%%%%%%%%%%%%%%%%%%%%%%%%%%%%%%%%%%%%%%%%%%%%%%%%%%%%%%

%%%%%%%%%%%%%%%%%%%%%%%%%%%%%%%%%%%%%%%%%%%%%%%%%%%%%%%
%%%%%%%%%%%%%%%%%%%%%%%%%%%%%%%%%%%%%%%%%%%%%%%%%%%%%%%

%%%%%%%%%%%%%%%%%%%%%%%%%%%%%%%%%%%%%%%%%%%%%%%%%%%%%%%%%
%%%%%%%%%%%%%%%%%%%%%%%%%%%%%%%%%%%%%%%%%%%%%%%%%%%%%%%%%
\end{document}